\documentclass[12pt,reqno]{amsart}
\usepackage{amsmath}
\usepackage{amsxtra}
\usepackage{amscd}
\usepackage{amsthm}
\usepackage{amsfonts}
\usepackage{amssymb}
\usepackage{eucal}

\newtheorem{theorem}{Theorem}[section]

\newtheorem{lem}[theorem]{Lemma}

\theoremstyle{definition}

\theoremstyle{remark}

\theoremstyle{remark}

\numberwithin{equation}{section}

\newcommand{\nc}{\newcommand}
\nc{\on}{\operatorname}
\nc{\ch}{\mbox{ch}}
\nc{\Z}{{\mathbb Z}}
\nc{\C}{{\mathbb C}}
\nc{\pone}{{\mathbb C}{\mathbb P}^1}
\nc{\pa}{\partial}
\nc{\F}{{\mathcal F}}
\nc{\arr}{\rightarrow}
\nc{\larr}{\longrightarrow}
\nc{\al}{\alpha}
\nc{\ri}{\rangle}
\nc{\lef}{\langle}
\nc{\W}{{\mathcal W}}
\nc{\la}{\lambda}
\nc{\ep}{\epsilon}

\nc{\su}{\widehat{{\mathfrak sl}}_2}
\nc{\sw}{{\mathfrak s}{\mathfrak l}}

\nc{\g}{{\mathfrak g}}
\nc{\h}{{\mathfrak h}}
\nc{\n}{{\mathfrak n}}
\nc{\N}{\widehat{\n}}
\nc{\G}{\widehat{\g}}
\nc{\De}{\Delta_+}
\nc{\gt}{\widetilde{\g}}
\nc{\Ga}{\Gamma}
\nc{\one}{{\mathbf 1}}
\nc{\z}{{\mathfrak Z}}
\nc{\zz}{{\mathcal Z}}
\nc{\Hh}{{\mathcal H}_\beta}
\nc{\qp}{q^{\frac{k}{2}}}
\nc{\qm}{q^{-\frac{k}{2}}}
\nc{\La}{\Lambda}
\nc{\wt}{\widetilde}
\nc{\qn}{\frac{[m]_q^2}{[2m]_q}}
\nc{\cri}{_{\on{cr}}}
\nc{\kk}{h^\vee}
\nc{\sun}{\widehat{\sw}_N}
\nc{\hh}{\widehat{\mathfrak h}}
\nc{\HH}{{\mathcal H}_{q,t}}
\nc{\ca}{\wt{{\mathcal A}}_{h,k}(\sw_2)}
\nc{\gl}{\widehat{{\mathfrak g}{\mathfrak l}}_2}
\nc{\el}{\ell}
\nc{\s}{{\mathbf s}}
\nc{\bi}{\bibitem}
\nc{\om}{\omega}
\nc{\WW}{\W_\beta}
\nc{\scr}{{\mathbf S}}
\nc{\ab}{{\mathbf a}}
\nc{\rr}{r}
\nc{\ol}{\overline}
\nc{\con}{qt^{-1} + q^{-1}t}
\nc{\den}{q^{\el-1} t^{-\el+1}+ q^{-\el+1} t^{\el-1}}
\nc{\ds}{\displaystyle}
\nc{\B}{B}
\nc{\A}{{\mathbb A}}
\nc{\GG}{{\mathcal G}}
\nc{\UU}{{\mathcal U}}
\nc{\MM}{{\mathcal M}}
\nc{\CC}{{\mathcal C}}
\nc{\GL}{{}^L G}
\nc{\dzz}{\frac{dz}{z}}
\nc{\Res}{\on{Res}}
\nc{\rep}{{\mathcal R}ep \;}
\nc{\uqg}{U_q \G}
\nc{\uqgg}{U_q \g}
\nc{\Fq}{{\mathbb F}_q}

\nc{\stimes}{\ltimes}
\nc{\K}{\hat{\mathcal K}}
\nc{\Ql}{\ol{\mathbb Q}_\ell}

\nc{\ga}{\gamma}
\nc{\PL}{{}^L P}
\nc{\E}{\mc E}
\nc{\mc}{\mathcal}
\nc{\mbf}{\mathbf}
\nc{\bb}{{\mathfrak b}}
\nc{\OO}{{\mc O}}
\nc{\Po}{{\mc P}}
\nc{\V}{{\mc V}}
\nc{\yy}{{\mc Y}}
\nc{\M}{\mathcal M}
\nc{\Coh}{{{\mathcal C}oh}}
\nc{\Cohn}{\Coh_n}
\nc{\f}{{\mathcal F}}
\nc{\si}{_E}
\nc{\Gaf}{{\mathbb G}_{a,\Fq}}
\nc{\KK}{{\mathfrak k}}

\nc{\PCr}{{ \bs P  (\C[x])^r   }}
\nc{\PCN}{{ (\bs P  \C[x])^N   }}
 
\nc{\sN}{sl_{2N+1}}

\nc{\Pzr}{{ \bs P(\C((x-z)))^r}}

\nc{\PzN}{{ \bs P(\C((x-z)))^N}}

\newcommand{\bean}{\begin{eqnarray}}
\newcommand{\eean}{\end{eqnarray}}
\newcommand{\be}{\begin{displaymath}}
\newcommand{\ee}{\end{displaymath}}
\newcommand{\bea}{\begin{eqnarray*}}   
\newcommand{\eea}{\end{eqnarray*}}
\newcommand{\bs}{\boldsymbol}
\newcommand{\Ref}[1]{{$($\ref{#1}$)$}}

\begin{document}

\title[Populations of Bethe solutions]
{Populations of solutions of the $XXX$ Bethe equations associated 
to Kac-Moody algebras}
\author[E. Mukhin and A. Varchenko]
{E. Mukhin and A. Varchenko}
\thanks{Research of E.M. is supported in part by NSF grant DMS-0140460.
Research of A.V. is supported in part by NSF grant DMS-9801582}
\address{E.M.: Department of Mathematical Sciences, Indiana University -
Purdue University Indianapolis, 402 North Blackford St, Indianapolis,
IN 46202-3216, USA, \newline mukhin@math.iupui.edu}
\address{A.V.: Department of Mathematics, University of North Carolina 
at Chapel Hill, Chapel Hill, NC 27599-3250, USA, anv@email.unc.edu}
\maketitle

\begin{abstract}
We consider the $XXX$ Bethe equation associated 
with integral dominant weights of
a Kac-Moody algebra and introduce a generating procedure
constructing new solutions starting from a given one. 
The family of all solutions constructed from a given one is
called a population. We describe properties of populations.
\end{abstract}

\section{Introduction}

The $XXX$ Bethe equation is a system of algebraic equations
associated to a Kac-Moody algebra $\g$, a non-zero step $h \in \C$,
distinct complex numbers $z_1,\dots,z_n$, integral dominant
$\g$-weights $\La_1,\dots,\La_n$ and a $\g$-weight $\La_\infty$, 
see \cite{OW} and Section \ref{Bethe eqn sec} below.  

In the theory of the Bethe ansatz one associates a Bethe equation to
an integrable model. Then each solution of the Bethe equation gives an
eigenvector of commuting Hamiltonians of the model.  The general
conjecture is that the constructed vectors form a basis in the space of 
states of the model. There is vast literature on this subject, see e.g. 
\cite{BIK, Fa, FT}.

The first step to the above mentioned conjecture is
to count the number of solutions of the Bethe equation.
One can expect that the number of solutions is equal to the dimension
of the space of states of the model.

In an attempt to approach that counting problem, in this note we relate to
each solution of the Bethe equation an object called the
population. We hope that it might be easier to count populations than
individual solutions. For instance, if the Kac-Moody algebra is
$sl_{r+1}$, then each population corresponds to a point of intersection
of suitable Schubert varieties in a suitable Grassmannian variety. Then 
the Schubert calculus allows us to count the populations, see
\cite{MV1}, \cite{MV2}.

The populations in the case of $sl_{r+1}$ were introduced in \cite{MV1}. 
The construction of a population goes as follows.  
First 
one defines $r$ polynomials $T_1(x), \dots , T_r(x)$, \
$T_i(x) = \prod_{s=1}^n\prod_{p=1}^{(\Lambda_s,\alpha_i)}
(x-z_s - (\La_1, \alpha_i)h + 2ph)$,
where $\alpha_1, \dots , \alpha_r$ are simple roots and
$(\ ,\ )$ is the standard scalar product.
A solution of the Bethe equation is an $r$-tuple of polynomials 
$\bs y = (y_1(x),\dots,y_r(x))$, each considered up to multiplication
by a non-zero number. The fact that $\bs y$ is a solution is
formulated as follows. The tuple $\bs y$ is a solution, if
for every $i=1, \dots, r$, there exists a polynomial
$\tilde y_i(x)$ satisfying the equation
\be
\tilde y_i(x+2h) y_i(x)-\tilde y_i(x) y_i(x+2h)= 
T_i(x) y_{i-1}(x+h)y_{i+1}(x+h) 
\ee
(we also need some requirements on $\bs y$ to be generic). 
It turns out that for all complex numbers $c$ (with finitely many
exceptions) the tuples 
$(y_1(x),\dots, \tilde y_i(x) + cy_i(x), \dots, y_r(x))$
are also solutions of the Bethe equation.
The collection of all $r$-tuples 
which can be constructed  starting from $\bs y$
by iteration of the above procedure 
is called the population of solutions originated at $\bs y$.

In this paper we suggest a similar generation procedure for an
arbitrary solution of a Bethe equation of 
a non-homogeneous $XXX$ model
associated with an arbitrary Kac-Moody algebra. The new feature of
such a generation procedure is that the shift in the $i$-th difference
equation depends on the length of the $i$-th simple root of the
Kac-Moody algebra.

A population of solutions is an interesting object. 
It is an algebraic variety. 
It is a finite-dimensional irreducible
algebraic variety, if the Weyl group of
the Kac-Moody algebra is finite. 
In this short article we present only elementary properties
of populations, see Section \ref{prop sec}. 
The main open questions are :
\begin{enumerate}
\item[$\bullet$] Describe a population as 
an algebraic
variety,
\item[$\bullet$] Compute the number of populations originated at solutions of
a given Bethe equation.
\end{enumerate} 

The populations related to the Gaudin model of a Kac-Moody algebra
 $\g$ were introduced in \cite{MV2}. In \cite{MV3}, see also \cite{F}, 
 it is proved that in the case of semi-simple Lie algebras
 every $\g$-population is isomorphic to the flag variety of the
 Langlands dual Lie algebra $\g^L$. It is conjectured in \cite{MV1} 
 that the same is true for all Kac-Moody algebras. 
 It is plausible that
 every $\g$-population introduced in this paper is isomorphic to
 the flag variety of the Langlands dual Kac-Moody algebra $\g^L$. In
 \cite{MV1},\cite{MV4} this conjecture is proved for 
Lie algebras of types A,D,E.

The reproduction procedure for the Bethe ansatz type equations was studied in 
some previous works, see \cite{MV1}, \cite{MV4}. Note that the equations 
of the B, C type in
\cite{MV1} differ from the B, C type Bethe equations of the present paper, and
the equations associated to a Kac-Moody algebra in \cite{MV4}
differ from the Bethe equations of the present paper if and only if there is 
an non-zero non-diagonal entry of the Cartan matrix different from $-1$.

We remark that all the constructions and statements of the present paper can 
be repeated for the $XXY$ Bethe equation.

\bigskip 

We thank V. Tarasov for pointing out the references related to the Bethe
ansatz equations and interesting discussions.
 
\section{Kac-Moody algebras}\label{Kac-Moody sec}

Let $A=(a_{ij})_{i,j=1}^r$ be a generalized  Cartan matrix, 
$a_{ii}=2$,
$a_{ij}=0$ if and only $a_{ji}=0$,
 $a_{ij}\in \Z_{\leq 0}$ if $i\ne j$. 
We  assume that $A$ is symmetrizable,  
there is a diagonal matrix $D=\on{diag}\{d_1,\dots,d_r\}$ 
with relatively prime positive integers $d_i$ such that $B=DA$
is symmetric.

Let $\g=\g(A)$ be the corresponding complex Kac-Moody 
Lie algebra (see \cite{K}, \S 1.2), 
$\h \subset \g$  the Cartan subalgebra.
The associated scalar product is non-degenerate on $\h^*$ and 
 $\on{dim}\h = r + 2d$,  where $d$ is the dimension of 
the kernel of the Cartan matrix $A$.

Let $\al_i\in \h^*$, $\al_i^\vee\in \h$, $i = 1, \dots , r$, be the
sets of simple roots, 
coroots, respectively. We have
\bea
d_i=(\al_i,\al_i)/2, \qquad
 a_{ij}=2(\al_i,\al_j)/(\al_i,\al_i)=\langle \al_j,\al_i^\vee\rangle, \qquad
\langle\la ,\al^\vee_i\rangle=2(\la,\al_i)/{(\al_i,\al_i)},
\eea
where $\la\in\h^*$. 
Let $\mathcal P = \{ \lambda \in \h^* \, |\, \langle\la
,\al^\vee_i\rangle \in \Z\}$. Let 
$\mathcal P^+ = \{ \lambda \in \h^* \, |\, \langle\la
,\al^\vee_i\rangle \in \Z_{\geq 0}\}$. Elements of $\mc P^+$ 
are called dominant integral weights.

The Weyl group $\mathcal W\in\on{End (\h^*)}$ is generated by 
reflections $s_i$, $i=1,\dots,r$, 
\be
s_i(\la)=\la-\langle\la,\al_i^\vee\rangle\al_i, \qquad \la\in\h^*.
\ee

Fix $\rho\in\h^*$ such that $\langle\rho,\al_i^\vee\rangle=1$,
$i=1,\dots,r$. We have $(\rho,\al_i)= (\al_i,\al_i)/2$.
We use the notation $w\cdot\la$ for the shifted action of the Weyl
group given by
\bea\label{shifted}
w\cdot\la=w(\la+\rho)-\rho,\qquad w\in \mathcal W,\;\la\in\h^*.
\eea

\section{The Bethe equations}\label{Bethe eqn sec}
Let $\bs\La=(\La_i)_{i=1}^n\in(\mc P^+)^n$,
$\bs z=(z_i)_{i=1}^n\in\C^n$,  $\bs
l=(l_i)_{i=1}^r\in\Z^r_{\geq 0}$,
$\bs  t = ( t_j^{(i)})^{j = 1,\dots ,l_i}_{i = 1,\dots,r}$. 

Set  $\bar {\bs \La}=(\La_1,\dots,\La_n,\La_\infty)$, where
\bean\label{wt at inf}
\La_\infty = \sum_{i=1}^n \La_i  -  \sum_{i=1}^r
l_i\al_i\in \mc P.
\eean

Fix a non-zero complex number $h$.
The $XXX$ Bethe equation associated with $\bs z, \bar{\bs \La}$ is the
following system of algebraic equations for variables $\bs t$, see
\cite{OW}:
\bean\label{Bethe}
\prod_{s=1}^n\frac{t_j^{(i)}-z_s+(\La_s,\al_i) h}
{t_j^{(i)}-z_s-(\La_s,\al_i) h}\,
\prod_{m=1}^r\prod_{k=1}^{l_m}
\frac{t_j^{(i)}-t_k^{(m)}-(\al_m,\al_i) h}
{t_j^{(i)}-t_k^{(m)}+(\al_m,\al_i) h}=-1,
\eean
where $i=1,\dots,r$, $j=1,\dots,l_i$.

Introduce 
\be
h_i=(\al_i,\al_i)h.
\ee
Then the Bethe equation \Ref{Bethe} can be written in the form
\be
\prod_{s=1}^n\frac{t_j^{(i)}-z_s+\langle\La_s,\al_i^\vee\rangle h_i/2}
{t_j^{(i)}-z_s-\langle\La_s,\al_i^\vee\rangle h_i/2}\,
\prod_{m=1,\dots,r,\atop m\neq i}\prod_{k=1}^{l_m}
\frac{t_j^{(i)}-t_k^{(m)}-a_{im} h_i/2}
{t_j^{(i)}-t_k^{(m)}+a_{im}h_i/2}\prod_{k=1,\dots,l_i, \atop k\neq j}
\frac{t_j^{(i)}-t_k^{(i)}-h_i}{t_j^{(i)}-t_k^{(i)}+h_i}=1.
\ee

The product of symmetric groups
$S_{\bs l}=S_{l_1}\times \dots \times S_{l_N}$ acts on the  set of
solutions of \Ref{Bethe}
permuting the coordinates with the same upper index. 
An $S_{\bs l}$ orbit of solutions of \Ref{Bethe} such that
\begin{itemize}
\item $t_j^{(i)}\neq t_k^{(i)}+ph_i$ for all $i,j,k$, $j\neq k$ and $p=0,1$,
\item $t_j^{(i)}\neq t_k^{(m)}-a_{im}h_i/2-ph_i$ for all $i,j,k,m$,
  $i\neq m$ and $p=1,\dots, -a_{im}$,
\item $t_j^{(i)}\neq z_s+\langle \La_s,\al_i^\vee\rangle h_i/2-ph_i$ for 
 all $i,j,s$  and $p=1,\dots,\langle \La_s,\al_i^\vee\rangle$,
\end{itemize}
is called a {\it Bethe solution associated to $\bs z, \bar{\bs \La}$}.

Note that not all solutions of equation \Ref{Bethe} are Bethe solutions.

In the case of $sl_{r+1}$, the Bethe solutions $\bs t$ are in
one-to-one correspondence with $2h$-critical points $\bar {\bs t}$
defined in Section 3.1 of \cite{MV1}. This correspondence is given by 
$\bar t_k^{(i)}=t_k^{(i)}-ih$, $b_s^{(i)}=i/2-(\La_s,\al_i)/2$ in
notation of \cite{MV1}. Note
that in the case of $B_r$ and $C_r$, the equations of Sections 6.1 and
7.1 in \cite{MV1} are different from \Ref{Bethe}.

\section{Reproduction of the Bethe solutions}\label{repr sec}

Given  $\bs t = (t^{(i)}_j)$, we form an $r$-tuple of polynomials
of a variable $x$, $\bs y=(y_1,\dots,y_r)$,
$y_i(x)=\prod_{j=1}^{l_i}(x-t_j^{(i)})$. We say that $\bs y$ {\it
  represents  $\bs t$}. 

We are interested only
in the zeros of $y_i$, and we consider those polynomials up to
multiplication by a non-zero complex number, so the $r$-tuple $\bs y$
is an element of $(\bf P\C[x])^r$, where $\bf P\C[x]$ is the
projectivization of the space of polynomials in $x$.

For $i=1,\dots,r$, introduce polynomials $T_i(x)$: 
\be
T_i(x)=\prod_{s=1}^n\prod_{p=1}^{\langle\La_s,\al_i^\vee\rangle}
(x-z_s-\langle\La_s,\al_i^\vee\rangle h_i/2+ph_i).
\ee

An $r$-tuple of polynomials $\bs y$ is called {\it generic} if $y_i(x)$ has
no multiple roots and no common roots with polynomials $y_i(x+h_i), T_i(x),
y_m(x+a_{im}h_i/2+ph_i)$ for $m=1,\dots,r$ and $p=1,\dots,-a_{im}$.
If $\bs y$ represents a Bethe solution, then $\bs y$ is generic.

\begin{theorem}\label{fer teor} 
Let $\bs y$ be a generic $r$-tuple of polynomials. The
  following three properties are equivalent:

(i) $\bs y$ represents a Bethe solution associated to $\bs z,\bar{\bs
  \La}$,

(ii) for $i=1,\dots, r$, the polynomial
\bean\label{div pol}
\left(\prod_{s=1}^n(x-z_s-\langle \La_s,\al_i^\vee \rangle h_i/2)\prod_{m,\
  a_{im}<0} y_m(x+a_{im}h_i/2)\right) y_i(x+h_i)+\notag\\
\left(\prod_{s=1}^n(x-z_s+\langle \La_s,\al_i^\vee \rangle h_i/2)\prod_{m,\
  a_{im}<0} y_m(x-a_{im}h_i/2)\right) y_i(x-h_i)
\eean
 is divisible by $y_i(x)$,

(iii) there exist polynomials $\tilde y_1,\dots,\tilde y_r$ such that
for $i=1,\dots r$,
\bean\label{wr eq}
y_i(x+h_i)\tilde y_i(x)-y_i(x)\tilde
y_i(x+h_i)=T_i(x)\prod_{m=1}^r\prod_{p=1}^{-a_{im}}
  y_m(x+a_{im}h_i/2+ph_i/2).
\eean
\end{theorem}
Note that the product over $p$ in the part (iii) equals to 1 in the
case of $a_{im}=0$ and in the case of $i=m$.
\begin{proof}
Since $y_i$ has no multiple roots, condition (ii) is equivalent to vanishing of the polynomial
\Ref{div pol} at the roots of $y_i$. This condition is exactly equation \Ref{Bethe}. Thus
(i) and (ii) are equivalent, cf. Lemma 2.2 in \cite{MV1}.

Consider the following equation for a function $c_i(x)$:
\be
c_i(x+h_i)-c_i(x)=\frac{T_i(x)\prod_{m=1}^r\prod_{p=1}^{-a_{im}}
  y_m(x+a_{im}h_i/2+ph_i/2)}{y_i(x)y_i(x+h)}.
\ee
Condition (iii) is equivalent to the existence of 
the solution of the form $c(x)=\tilde y_i(x)/y(x)$, where $\tilde y_i(x)$ is a polynomial. 
This is equivalent to the condition on residues at zeroes of $y_i$:
\be
\left(\on{Res}_{x=t_j^{(i)}}+\on{Res}_{x=t_j^{(i)}-h_i}\right)\left(\frac{T_i(x)
\prod_{m=1}^r\prod_{p=1}^{-a_{im}}
y_m(x+a_{im}h_i/2+ph_i/2)}{y_i(x)y_i(x+h)}\right)=0,
\ee
which coincides with equation \Ref{Bethe}. Thus (i) and (iii) are equivalent, cf. Lemmas 2.3 
and 2.4 in \cite{MV1}.
\end{proof}

An $r$-tuple $\bs y$ 
satisfying condition (iii) of the theorem is called {\it
  fertile with respect to $\bs z, \bs \La$}. 
An $r$-tuple $\bs y$ represents a Bethe solution if and only if it
is fertile and generic. If $\bs y$ is fertile, then the $r$-tuples of the form 
$\bs y^{(i)}=(y_1,\dots,\tilde y_i,\dots,y_r)$ are called the {\it
  immediate descendants of $\bs y$ in the $i$-th direction}. 

Let $\C_d[x]$ be the space of all polynomials of degree at most $d$.
The set of fertile tuples is closed in $(\bs P \C_d[x])^r$:

\begin{lem}\label{limit fertile lem}
Assume that a sequence  of fertile tuples of 
polynomials $\bs y_k$, $k = 1, 2, \dots$, 
has a limit $\bs y_\infty$ in $(\bs P \C_d[x])^r$ as $k$
tends to infinity. 
Then the limiting tuple $\bs y_\infty$ is fertile.

Assume in addition, that for some $i$,
all immediate descendents of all $y_k$ in the $i$-th direction are in 
$(\bs P \C_d[x])^r$.
If $\bs y^{(i)}_\infty$ is 
an immediate descendant of $\bs y_\infty$ 
in the $i$-th direction, then $\bs y^{(i)}_\infty\in (\bs P \C_d[x])^r$ and
there exist immediate descendants
$\bs y_k^{(i)}$ of $\bs y_k$ in the $i$-th direction such that 
$\bs y^{(i)}_\infty$ is the limit of $\bs y_k^{(i)}$.

\end{lem}

\begin{proof} (Cf. Lemma 3.4 in \cite{MV1}.) 
For each $k$ we have a plane of polynomials $P_k^{(i)}$ spanned 
by $(\bs y_k)_i$ and $({\bs y}_k^{(i)})_i$. Since the Grassmannian variety of planes 
in $\C_d[x]$ is compact, there is a limiting plane $P_\infty^{(i)}$. This plane $P_\infty^{(i)}$ 
obviously contains 
$(\bs y_\infty)_i$ and the lemma follows.
\end{proof}

Now we are ready to 
prove
the main result of this section.

\begin{theorem}\label{simple rep thm}
Let $\bs y$ represent a Bethe solution associated to $\bs z, \bar{\bs
  \La}$, let $i\in\{1,\dots,r\}$ and let an 
 $r$-tuple $\bs y^{(i)}=(y_1,\dots,\tilde
  y_i,\dots,y_r)$ be an immediate descendant in the direction $i$.
  Assume that $\bs y^{(i)}$ is generic. Then 
$\bs  y^{(i)}$ represents a Bethe solution associated either to 
$\bs z, \bar{\bs \La}$ if $\deg \tilde y_i=\deg y_i$ or  to 
$\bs z, (\La_1,\dots,\La_n,s_i\cdot \La_\infty)$ if 
$\deg \tilde y_i\neq\deg y_i$.
\end{theorem}
\begin{proof}
Denote $\tilde t_j^{(i)}$ the zeroes of $\tilde y_i$.
We have to check equation \Ref{Bethe} for $\bs y^{(i)}$. 

Equation \Ref{Bethe} with respect to $t_j^{(k)}$ where $a_{ik}=0$ is the same 
for $\bs y$ and $\bs y^{(i)}$ and therefore it is satisfied.

Equation \Ref{Bethe} with respect to $\tilde t_j^{(i)}$ is equivalent to the 
existence of
a polynomial $\tilde{\tilde y}_i$ satisfying the $i$-th Wronskian 
identity as in \Ref{wr eq}, see Theorem \ref{fer teor}. We can simply take
$\tilde{\tilde y}_i=-y_i$.

Let now $k$ be such that $a_{ik}<0$.
To obtain the Bethe equation \Ref{Bethe} with respect to $t_j^{(k)}$ 
for $\bs y^{(i)}$ from the Bethe equation
for $\bs y$ with respect to the same variable $t_j^{(k)}$  it is enough to show that  
\be
\prod_{s=1}^{l_i}
\frac{t_j^{(k)}-t_s^{(i)}-a_{ki} h_k/2}
{t_j^{(k)}-t_s^{(i)}+a_{ki} h_k/2}=\prod_{s=1}^{{\tilde l}_i}
\frac{t_j^{(k)}-{\tilde t}_s^{(i)}-a_{ki} h_k/2}
{t_j^{(k)}-{\tilde t}_s^{(i)}+a_{ki}h_k/2}.
\ee
This condition can be rewritten as:
\bean\label{need}
\frac{y_i(t_j^{(k)}-a_{ki}h_k/2)}{y_i(t_j^{(k)}+a_{ki}h_k/2)}=\frac{\tilde y_i(t_j^{(k)}-a_{ki}h_k/2)}{\tilde y_i(t_j^{(k)}+a_{ki}h_k/2)}.
\eean
By substituting $x=t_j^{(k)}-a_{ik}h_i/2-ph_i/2$ where $p=1,\dots,-a_{ik}$ 
in the $i$-th Wronskian identity 
\Ref{wr eq}, we observe that the right hand side vanishes and therefore we obtain
\bean\label{have}
\frac{y_i(t_j^{(k)}-a_{ki}h_k/2-(p-1)h_i/2)}{y_i(t_j^{(k)}-a_{ki}h_k/2-ph_i/2)}
=\frac{\tilde y_i(t_j^{(k)}-a_{ki}h_k/2-(p-1)h_i/2)}{\tilde y_i(t_j^{(k)}-a_{ki}h_k/2-ph_i/2)}.
\eean
The identity \Ref{need} is obtained by the product of identities \Ref{have} for 
$p=1,\dots,-a_{ik}$, since $a_{ik}h_i=a_{ki}h_k=2(\al_i,\al_k)h$.

Finally, if $\deg \tilde y_i\neq\deg y_i$ then using the $i$-th Wronskian
identity \Ref{wr eq} once again we obtain a relation
\be
\deg y_i+\deg \tilde y_i=1+\deg T_i-\sum_{m, a_{im}<0}a_{im}\deg y_m,
\ee
which implies the relation between the weights at infinity stated in the theorem.
\end{proof}

Note that if the polynomial $\tilde y_i$ satisfies \Ref{wr eq} then
for any complex number $c$, the polynomial $\tilde y_i+cy_i$ also
satisfies \Ref{wr eq}. For
all but finitely many values of $c$ the $r$-tuple $(y_1,\dots, \tilde
y_i+cy_i,\dots, y_r)$ is generic. Thus starting from a Bethe solution
we constructed a family of $r$-tuples of polynomials, all but finitely
many of which represent Bethe solutions. Then by Lemma \ref{limit fertile
  lem}, all these $r$-tuples are fertile.

We call this construction the {\it simple reproduction in the
  $i$-th direction}.

\section{Populations of Bethe solutions}\label{prop sec}
Let $\bs y_0$ represent a Bethe solution associated to $\bs z,
\bar{\bs \La}$. 

We apply the simple reproduction
procedure in each direction $i$ and 
obtain a family of fertile $r$-tuples $\bs y_0^{(i)}$. Then we
  apply the simple reproduction procedure in each direction $j$ 
 to these tuples again and obtain
  a larger family of $r$-tuples $(\bs y_0^{(i)})^{(j)}$. We continue
  this procedure which we call the {\it reproduction procedure}. The
  set of all $r$-tuples obtained by the reproduction procedure from $\bs
  y_0$ is called the {\it populations of Bethe solutions originated at
  $\bs y_0$} and is denoted $P(\bs y_0)$.

More formally, the population of $\bs y_0$ is the set of all $r$-tuples
of polynomials ${\bs y}\in ({\bf P}\C[x])^r$ such that there
exist $r$-tuples $\bs
y_1$, $\bs y_2,\dots , \bs y_k={\bs y}$, where $\bs y_i$
is an immediate descendant of $\bs y_{i-1}$ for all $i=1,\dots,k$.

The following lemma is straightforward.
\begin{lem}
All $r$-tuples in the population $P(\bs y_0)$ are fertile with respect to $\bs z,\bs
\La$. \hfill $\Box$
\end{lem}

For an $r$-tuple ${\bs y}$, define the corresponding 
 weight at infinity 
by formula \Ref{wt at inf} where we
  set $l_i=\deg {y}_i$. From Theorem \ref{simple rep thm}, we
  obtain
\begin{lem}\label{weights at inf}
The set of weights at infinity of the $r$-tuples in the population $P(\bs
y_0)$ coincides with the orbit of $\La_\infty$ under the shifted action
of the Weyl group of $\g$.\hfill $\Box$
\end{lem}

In particular, if ${\bs y}\in P(\bs y_0)$ is generic then it
represents a Bethe solution associated to $\bs z,
(\La_1,\dots,\La_n,w\cdot \La_\infty)$ for some element $w$ of the
Weyl group $\mc W$.

If $\bs y$ represents a Bethe solution associated to $\bs z,\bar{\bs
  \La}$, then we say that the population $P(\bs y)$ {\it is associated to
  $\bs z, \bs \La$}.
Clearly, 

\begin{lem}\label{int=same} If two populations associated to $\bs z,\bs \La$ 
intersect then they coincide.  \hfill $\Box$
\end{lem}

Let $\bs y$ and $\bar {\bs y}$ represent Bethe solutions. 
We say that $\bs y\equiv \bar {\bs y}$ if and only if $\bs y$ and $\bar {\bs y}$ 
belong to the same population. Then by Lemma \ref{int=same}, $\equiv$ 
is an equivalence relation.

The following lemma follows from the general standard arguments, 
cf. Corollaries 3.13 and 3.14 in \cite{MV2}. 

\begin{lem}
For any $d\in\Z_{\geq 0}$, the intersection of any population 
of Bethe solutions with $(\bf
P\C_d[x])^r$ is an algebraic variety. In particular, if the Weyl group
of $\g$ is finite then the population of Bethe solutions is an irreducible algebraic
variety. \hfill $\Box$ 
\end{lem}

In conclusion we describe situations in which 
the number of populations can be computed by general arguments.

The $r$-tuple $(1, \dots , 1)\in (\bf P\C[x])^r$ is the unique $r$-tuple of
non-zero polynomials of degree 0. 
The weight at infinity of $(1, \dots , 1)$ is
$\Lambda_{\infty, (1, \dots , 1)} 
= \sum_{s=1}^n \La_s$. 
Let $P_{(1, \dots , 1)}$ be the population associated  to the initial data 
and originated at $(1, \dots , 1)$.

The next lemma says that if  $\Lambda_{\infty, (1, \dots , 1)}$ is in the orbit of $\La_\infty$  
then there is 
exactly one population of Bethe solutions, cf. Corollary 4.4 in
\cite{MV1} and Corollary 3.17 in \cite{MV2}.

\begin{lem}
Let $\bs y$ represent a Bethe solution such that
$\La_\infty$ has the form $w\cdot\Lambda_{\infty, (1, \dots , 1)}$ for some
$w\in\mc W$. Then $\bs y$ belongs to the population  
$P_{(1, \dots , 1)}$. 
\end{lem}
\begin{proof}
This lemma follows directly from Lemmas \ref{weights at inf} and \ref{int=same}.
\end{proof}

The last two lemmas give some sufficient conditions for the absence of Bethe
solutions, cf. Corollary 4.3 in \cite{MV1} and  
Corollaries 3.15, 3.16 in \cite{MV2}.

\begin{lem}\label{no1}
Suppose there is an element $w$ 
of the Weyl group such that $\sum_{s=1}^n \La_s\ -\ w \cdot  \La_\infty$ 
does not belong to the cone $\Z_{\geq 0} \al_1\ \oplus\ \dots\ \oplus \ 
\Z_{\geq 0} \al_r$. Then there are no Bethe solutions associated to
$\bs z,\bar{\bs\La}$. 
\end{lem} 
\begin{proof}
 The Lemma \ref{no1} follows from Lemma \ref{weights at inf}
 and the absence of polynomials of negative degrees.
\end{proof}

\begin{lem}\label{no2}
Suppose that $s_i\cdot \La_\infty=\La_\infty$ for some simple
reflection $s_i$. Then there are no Bethe solutions associated to
$\bs z,\bar{\bs\La}$.
\end{lem} 
\begin{proof}
If $\bs y$ represents a Bethe solution, then it is fertile in all directions. 
Then there is a polynomial $\tilde y_i$ satisfying \Ref{wr eq}
which has a degree different 
from the degree of $y_i$. Therefore the corresponding weights at infinity for 
$\bs y$ and $\bs y^{(i)}$ are different. This contradicts the
assumption of the 
lemma in view of Theorem \ref{simple rep thm}.
\end{proof}

\end{document}